\documentclass[a4paper,12pt,reqno,showkeys,superscriptaddress,nofootinbib]{revtex4}
\usepackage[centertags]{amsmath}
\usepackage{amsfonts}
\usepackage{amssymb}
\usepackage{amsthm}
\usepackage{newlfont}
\usepackage{stmaryrd}
\usepackage{mathrsfs}
\usepackage{mathtools}
\usepackage{euscript}
\usepackage{graphicx}
\usepackage{enumerate}
\usepackage{todonotes}

\usepackage{dblfloatfix}
\usepackage{color}

\usepackage{tikz}
\usepackage{pgf}
\usetikzlibrary{positioning,fit,calc}
\usetikzlibrary{arrows,automata}
\usepackage{wrapfig}
\usepackage{subfigure}
\usepackage{amscd}
\usepackage{hyperref}
\usepackage{cancel}

\theoremstyle{plain}
  \newtheorem{theorem}{Theorem}[section]
  
  \newtheorem{proposition}[theorem]{Proposition}
  \newtheorem{lemma}[theorem]{Lemma}
\theoremstyle{definition}

\theoremstyle{remark}

\usepackage{import}

\numberwithin{equation}{section}

\let\ve=\varepsilon

\newcommand{\bbL}{{\mathbb L}}

\newcommand{\bbV}{{\mathbb V}}

\newcommand{\bbZ}{{\mathbb Z}}

\newcommand{\opunit}{\text{1}\kern-0.22em\text{l}}

\DeclareMathAlphabet{\mathpzc}{OT1}{pzc}{m}{it}

\newcommand{\fig}{Fig.\;}

\newcommand{\id}{\textrm{d}}

\usepackage{floatrow}
\usepackage{caption}

\begin{document}
\title{Drazin-inverse and heat capacity for \\driven random walkers on the ring}
\author{Faezeh Khodabandehlou}
	\affiliation{Instituut voor Theoretische Fysica, KU Leuven, Belgium}
\author{Irene Maes}
	\affiliation{Departement Wiskunde, KU Leuven, Belgium}	
	\email{irene.maes@kuleuven.be}
\keywords{pseudo-inverse, quasipotential, nonequilibrium heat capacity.}

\begin{abstract}
We apply single and double tree-like representations of Markov jump processes on $\bbZ_N$ for obtaining their nonequilibrium heat capacity and for taking the diffusion limit $N\uparrow \infty$.  The main tool is a graphical representation of the Drazin-inverse of the backward generator. In that way, the combination of algebraic and graph-theoretical approaches enables an exact computation of an important nonequilibrium quantity (excess heat) for Markov jump processes.
\end{abstract}
\maketitle

\section{Introduction}
Markov jump processes offer an important modeling environment for nonequilibrium statistical mechanics.  They arise from quantum processes for multilevel systems (beyond decoherence times) or provide coarse-grained descriptions of diffusion processes, \cite{ja}. Nonequilibrium features such as maintaining steady currents, are possible by the presence of loops in the graph of possible configurations.  The present paper focuses on such a loop-dynamics and how one can extract the thermal response.\\

Our constructive approach to exact computations in nonequilibrium calorimetry combines elements of linear algebra (pseudo-inverses), graph theory (matrix-forest theorem) and Markov jump processes.  The purpose is to describe physically relevant properties of driven jump processes related to so called excess heat.  
In that way, the paper takes up new challenges both in statistical mechanics and within the context of stochastic processes; see also \cite{heatb,ner}. \\ 
More specifically, we are concerned with biased (or driven) random walks in an energy landscape defined on a ring.  Without the driving, the stationary process is reversible with respect to the usual Gibbs distribution for that energy function.  With the driving, we are motivated by nonequilibrium calorimetry  for computing the heat capacity of such driven systems, \cite{noncal}.  As we will see, such a computation depends on getting the stationary distribution and the excess work (the so-called quasipotential).  For both, we give algorithms following a graphical representation.  It allows fast computation and visualization of the relevant physics, in particular, to find the nonequilibrium heat capacities.\\

In the next section, we define the quasipotential in the context of Markov jump processes, and we explain the main questions.  In Section \ref{generator}, we show that the calculation of the quasipotential reduces to finding the Drazin-inverse of the backward generator.  That Drazin-inverse coincides here with the limit of the resolvent, which, from the matrix-forest theorem \cite{MForest1,MForest2}, gives rise to a graphical algorithm. That is explicitly achieved in Section \ref{ring} with the main theorem in Section \ref{grre}.  In that way, we achieve a fast calculation of quasipotentials in Section \ref{dop}; see also \cite{simon}.
Combined with the matrix-tree theorem \cite{kir}, it leads to an exact computation of the heat capacity in Section \ref{hc}.  The explicit low-temperature behavior may or may not show a Third Law of Thermodynamics (vanishing of the heat capacity with temperature).  Related to that, we discuss  boundedness of the quasipotential in Section \ref{nera}.  
Section \ref{stat} presents the Kirchhoff formula, which allows a diffusion limit to recover the stationary distribution for the corresponding driven diffusion on the circle,  as was already obtained in \cite{davide}. The same is achieved for the quasipotential in Section \ref{dp}, and in Section \ref{nera} we show how to use the graphical representation for obtaining bounds on the quasipotential. Conclusions are drawn in Section \ref{con}.\\
Appendix \ref{alg} is devoted to an explanation of the ideas in the algorithm and the computer code (which is publicly available in \cite{code}).  In the Appendix \ref{lap} we give an explicit proof that the index of the backward generator equals one.

\section{Setup and main results}\label{setup}
Consider  a finite connected graph $\cal G(\cal V(\cal G), \cal E(\cal G))$, without self-edges, i.e., there is no edge that starts and ends in the  same vertex. $\cal V(\cal G)$ is the set of  vertices and  $\cal E(\cal G)$ is the set of edges; put  $n= |\cal V(\cal G)|$ and $m= |\cal E(\cal G)|$. To every oriented edge $(x,y)$ a weight $k(x,y)>0$ is assigned.\\
A Markov jump process $X_t$ uses the set of vertices as possible (physical) states and the weights on the graph become the transitions rates over the edges:  $k(x,y)$ is the rate of jumping from $x$ to $y$ and no jump is possible when no edge is present between the two vertices.\\
The backward generator $L$ for the Markov process is the $n\times n$ matrix with elements
\begin{equation}\label{L}
    L_{x,y}=\begin{cases}
    L_{x,y}=k(x,y)&\quad x\not =y \\
    L_{x,x}=-\sum_y k(x,y)
    \end{cases}
\end{equation}
for $x,y \in \cal V(\cal G)$.  $L$ is directly related to the Laplacian matrix of the weighted directed graph; see also Appendix \ref{lap} and the many references in \cite{hindustan}.  For a weighted graph the Laplacian matrix is also called the Kirchhoff matrix; see e.g. \cite{che2006}. \\
Every function $ g$ on $ \cal V(\cal G)$ is a column vector. The semigroup $S(t) = \exp (tL), t\geq 0$,  satisfies
\[
S(t)g(x) = \langle g(X_t)\,|\, X_0=x\rangle  
\]
where  $\langle\cdot \,|\, X_0=x\rangle$ is the process expectation of the random walk $X_t$ starting from $X_0=x$ for  $x\in \cal V(\cal G)$.\\
The Master equation is 
\[\frac{d}{dt}\rho_t=\rho_t L\]
for the time-dependent probability $\rho_t$, written as a row vector.  There is a unique stationary distribution $\rho^s>0$, solution of $\rho^s L=0$, and, to settle the notation,
\[
\lim_{t \to \infty} \langle g(X_t)\,|\, X_0=x\rangle = \langle g\rangle^s := \sum_{y\in \cal V(\cal G)}\rho^s(y) g(y)
\]
where the convergence is exponentially fast, uniformly in $x$.\\

Take a source function $f$  on $\cal V(\cal G)$ with zero stationary expectation $\langle f\rangle^s =0$, and define
the quasipotential
\begin{equation}\label{pp}
V_f = \int_0^\infty \id t\, S(t)f
\end{equation}
which is exponentially converging.  The main subject of the paper is to use an algorithm to find $V_f$.  Note that equation \eqref{pp} can be rewritten as  $LV_f=-f$, which is a Poisson equation.  In the literature on that Poisson equation, the requirement  $\langle f\rangle^s =0$ is called a {\it centrality condition}, and most applications are related to diffusion processes; see e.g. \cite{poiss}.  Here we deal with the Poisson equation in a discrete setting. \\

The following main results are obtained in the coming sections:\\
 (1) identifying the quasipotential $V_f =-L_D^{-1} f$ as obtainable from the Drazin-inverse $L_D^{-1}$ of $L$ working on $f, \langle f\rangle=0$;\\
 (2) a graphical representation of the quasipotential $V$ on the ring, and its efficient computation for different $\ve, N, T$;\\
 (3)  from the graphical algorithm: exact computation of the nonequilibrium heat capacity $C$ as function of temperature $T$ (including $T\downarrow 0$)  with parameters $N,\ve$ on the ring;\\
(4) diffusion limit of Kirchhoff formula and quasipotential on the ring.

\section{The quasipotential is a pseudo-inverse}\label{generator}
The mathematical challenges of Markov jump processes on finite connected graphs often relate to questions of linear algebra.\\ 
The backward generator in equation \eqref{L} is a singular matrix, which invites the study of its pseudo-inverses.  One such pseudo-inverse is the Drazin inverse.\\ 

Consider a square matrix $A \in \mathbb{C}^{n\times n}$.
The \textit{index} of $A$, $\text{ind}(A)$, is the smallest non-negative integer $k$ such that $\text{rank}(A^k) = \text{rank}(A^{k+1})$.\\
The Drazin inverse of  $A$ is the matrix $A_D^{-1}$ satisfying the following three conditions \cite{book}: 
\begin{enumerate}\label{Drazin}
	\item $A^k A_D^{-1} A = A^k$ with $k=\text{ind}(A)$,
	\item $A_D^{-1} A A_D^{-1} = A_D^{-1}$,
	\item $A A_D^{-1} = A_D^{-1} A$.
\end{enumerate} 
The Drazin-inverse $A_D^{-1}$ always exists and is unique, \cite{book}.
We show here that the quasipotential \eqref{pp} can be obtained as the Drazin-inverse of the backward generator $L$. 

\begin{proposition}\label{p1}
	Let $V_f$ be the quasipotential \eqref{pp}  for a function $f$ having $\left\langle f\right\rangle^s = 0$, then 
\[
 V_f = -L_D^{-1} f.
\]
\end{proposition}
\begin{proof}
	In Appendix \ref{lap} we discuss why $\text{ind}(L) = 1$, and we give the explicit proof for the dynamics on a ring; see also \cite{che2006}. As follows from Theorem 3.1 (2) in \cite{Laurent}, that implies that
 \begin{equation*}
	\lim_{\lambda \to 0}(L + \lambda I)^{-1}L = L_D^{-1} \,L .
	\end{equation*}
 Because of the exponential convergence in \eqref{pp}, $\langle V_f\rangle^s = 0 $ and 
 \[
V_f(x) = \lim_{\lambda \to 0}\dfrac{-1}{(L - \lambda I)}f(x)
\]
where $\dfrac{1}{L-\lambda I}, \lambda \in \mathbb{C}, 0 < |\lambda|<\delta$ for some $\delta>0$, is the resolvent-inverse of the backward generator $L$.  In fact, $V_f$ is the unique solution of the equation $LV=-f$.  Therefore,
	\begin{align*}
	\lim_{\lambda \to 0} [(L + \lambda I)^{-1} f] &= -\lim_{\lambda \to 0} [(L + \lambda I)^{-1} LV]  \\
		& = -L_D^{-1} L V_f\\
	& = L_D^{-1} f.
	\end{align*}
 \end{proof}
 
As a remark, if the index of a matrix is $1$, then the special case of the Drazin-inverse is known as the group inverse. Therefore,  
 the Drazin--inverse $L_D^{-1}$ also equals the group inverse of $L$; see e.g. Theorem 2.2.1 in \cite{book}.

\section{Quasipotential on a ring}\label{ring}
Consider the ring $\frac 1{N}\bbZ_N \cong \bbZ_N$ with $N$ vertices ${0,\frac{1}{N},\frac{2}{N},\dots,\frac{N-1}{N}}$ . Moving in the clockwise direction takes $x$ to $x+1/N$ with driving $\ve$. 
The reason for considering  a mesh size $1/N$ is to be able to explore as well the diffusion limit where $N\uparrow \infty$ in what is to come.\\
A random walk $X_t$ on $\bbZ_N$ has transition rates $k\left(x,x+\frac{1}{N}\right)$ for the clockwise jumps $x\rightarrow x+\frac{1}{N}$ and has transition rates $k\left(x,x-\frac{1}{N}\right)$ for the counter-clockwise jumps $x\rightarrow x-\frac{1}{N}$; all other transitions are forbidden.\\

For an elementary example, we look at the simple random walker on the ring, where all the transition rates are equal to one and we take a source function $f$ on $\bbZ_N$ with  $\sum f(x) =0$.  Then, the quasipotential  $V=V_f$ in equation \eqref{pp} satisfies $\sum_x V(x) =0$ and
\begin{equation}\label{vv}
V(x+ 1/N) + V(x-1/N) - 2V(x) = -f(x)
\end{equation}
which is the standard discrete Poisson equation (for total charge equal to zero) with periodic boundary conditions, easily solved by Fourier transform.
To be explicit, for $N=3$, the solution is $V_f(x) = f(x)/3$. 

\subsection{Single and double trees}\label{mft}

\begin{figure}[h!]
    \centering
    \def\svgwidth{0.5\linewidth}        
    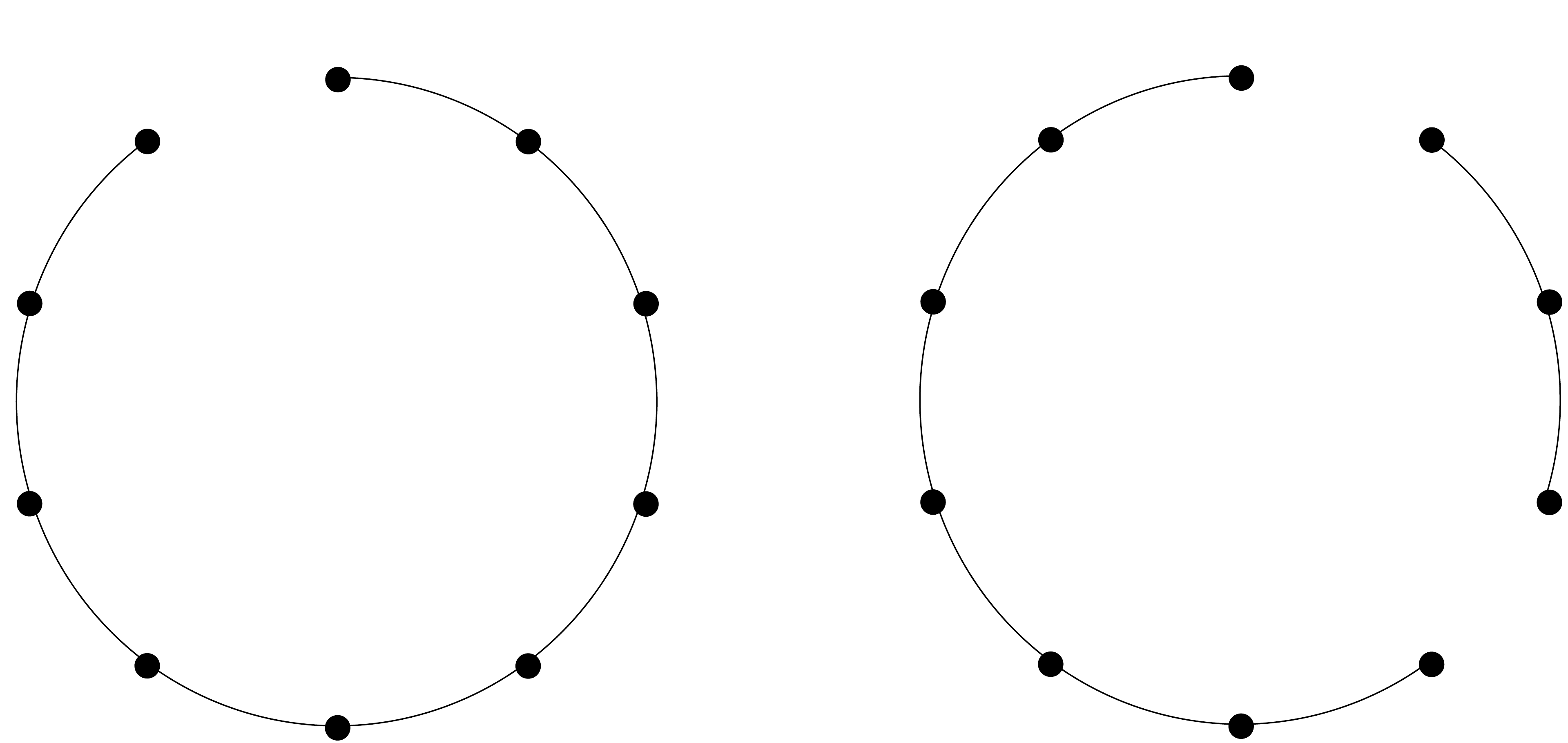
    \caption{{\small A tree in $\mathbb{Z}_N$ on the left side and double tree on the right side.}}\label{tree and directed tree}
\end{figure}

We will give a graphical representation of the quasipotential  \eqref{pp} for a general (nearest neighbor) random walk dynamics on the ring $\bbZ_N$.  We need to recall some definitions.\\
\textit{A tree} in the ring $\mathbb{Z}_N$ (see Fig.~\ref{tree and directed tree}) is the ring with exactly one edge removed. \textit{A double tree} in the ring $\mathbb{Z}_N$ (see Fig.~\ref{tree and directed tree}) is the ring with exactly two edges removed. Remark that a double tree consists of two connected parts and that one of the two parts is possibly only a vertex (if two edges with a common vertex are removed). 

\begin{figure}[h!]
    \centering
    \def\svgwidth{0.55\linewidth}        
\begingroup%
  \makeatletter%
  \providecommand\color[2][]{%
    \errmessage{(Inkscape) Color is used for the text in Inkscape, but the package 'color.sty' is not loaded}%
    \renewcommand\color[2][]{}%
  }%
  \providecommand\transparent[1]{%
    \errmessage{(Inkscape) Transparency is used (non-zero) for the text in Inkscape, but the package 'transparent.sty' is not loaded}%
    \renewcommand\transparent[1]{}%
  }%
  \providecommand\rotatebox[2]{#2}%
  \newcommand*\fsize{\dimexpr\f@size pt\relax}%
  \newcommand*\lineheight[1]{\fontsize{\fsize}{#1\fsize}\selectfont}%
  \ifx\svgwidth\undefined%
    \setlength{\unitlength}{1984.2519685bp}%
    \ifx\svgscale\undefined%
      \relax%
    \else%
      \setlength{\unitlength}{\unitlength * \real{\svgscale}}%
    \fi%
  \else%
    \setlength{\unitlength}{\svgwidth}%
  \fi%
  \global\let\svgwidth\undefined%
  \global\let\svgscale\undefined%
  \makeatother%
  \begin{picture}(1,0.45714286)%
    \lineheight{1}%
    \setlength\tabcolsep{0pt}%
    \put(0,0){\includegraphics[width=\unitlength,page=1]{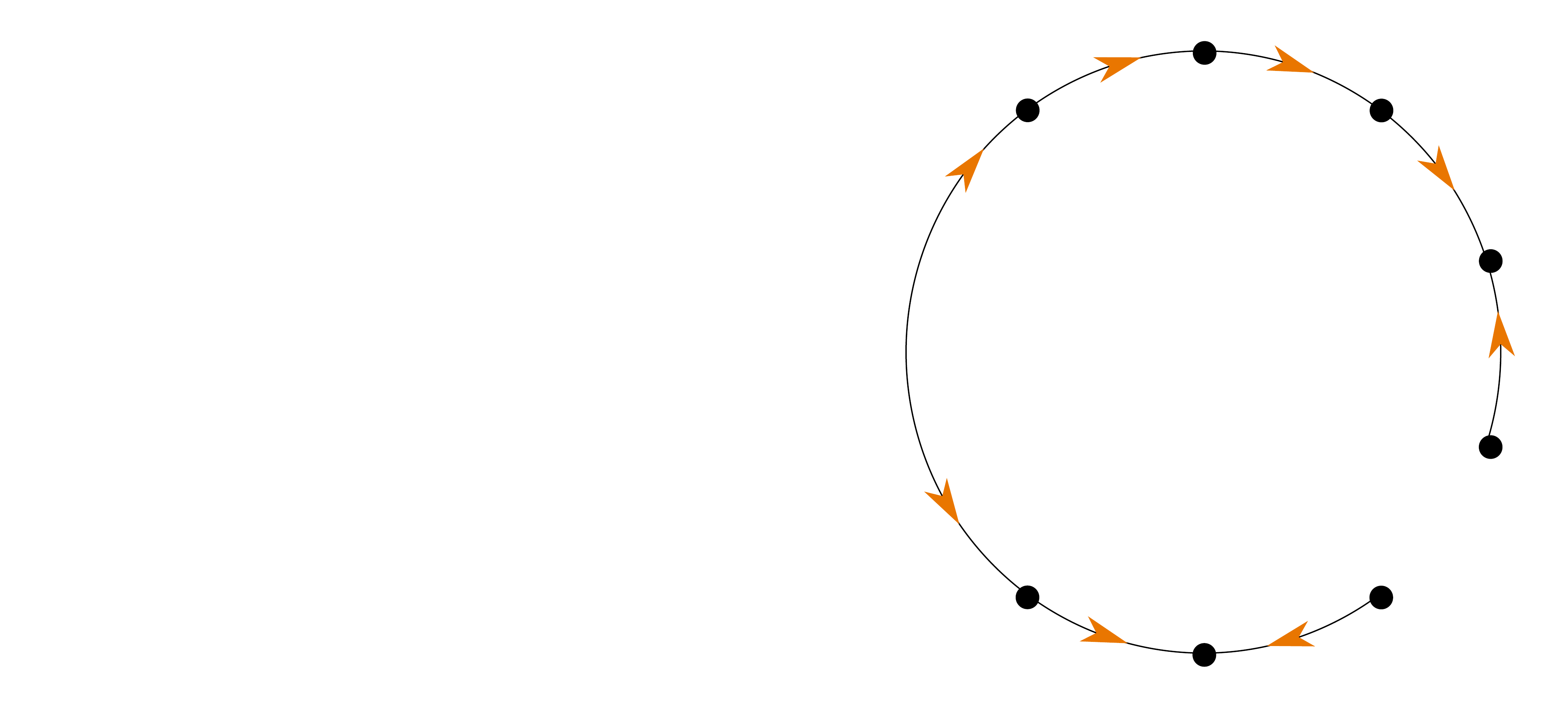}}%
    \put(0.63542389,0.40009421){\color[rgb]{0,0,0}\makebox(0,0)[lt]{\lineheight{1.25}\smash{\begin{tabular}[t]{l}$x$\end{tabular}}}}%
    \put(0.96844896,0.31024273){\color[rgb]{0,0,0}\makebox(0,0)[lt]{\begin{minipage}{0.08547835\unitlength}\raggedright $y$\end{minipage}}}%
    \put(0.75217765,0.0002056){\color[rgb]{0,0,0}\makebox(0,0)[lt]{\lineheight{1.25}\smash{\begin{tabular}[t]{l}$r$\end{tabular}}}}%
    \put(0,0){\includegraphics[width=\unitlength,page=2]{double_tree_and_directed_rooted_double_tree.pdf}}%
    \put(0.39816816,0.16640189){\color[rgb]{0,0,0}\makebox(0,0)[lt]{\lineheight{1.25}\smash{\begin{tabular}[t]{l}$r$\end{tabular}}}}%
    \put(0,0){\includegraphics[width=\unitlength,page=3]{double_tree_and_directed_rooted_double_tree.pdf}}%
  \end{picture}%
\endgroup%

    \caption{{\small A rooted tree with root in $r$ on the left side and a double rooted tree on the right side. The double rooted tree is an element of $\cal F^{x \to y}$.}}\label{double tree}
\end{figure}

\textit{A  rooted tree} in the ring (see Fig.~\ref{double tree}) is a tree in the ring where every edge has a direction and with a special vertex, called the root. The root is such that if we follow the directions from an arbitrary vertex on the ring, we end up in the vertex which is the root.

\textit{A double rooted tree} in the ring (see Fig.~\ref{double tree}) is a double tree in the ring where every edge has a direction and in each of the two parts of the double tree there is a special vertex called the root. This root is again such that if we start from an arbitrary vertex in one of the two parts and we follow the directions, we end up in the root of that part. 

Let $\cal T$ be the set of all rooted trees in the ring and $\cal T_y$ be the set of all rooted trees with root in $y$. Define $\cal F^{x \to y}$ to be the set of all double rooted trees in the ring such that $y$ is one of the roots and $x$ is in the same part as $y$ (so there is a path from $x$ to $y$, see Fig.~\ref{double tree}).\\ 
The weight of a rooted tree $T$ resp.~a double rooted tree $F$ in the ring, noted by $w(T)$ resp.~$w(F)$ is the product of all transition rates $k(z,u)$ such that the edge $\left(z,u\right)$ with direction $z \to u$ belongs to the tree or the double tree ($u$ is either $z+1/N$ or $z-1/N$).
If $S$ is some set of rooted trees or double rooted trees, we define $w(S)$ to be the sum of all weights of the elements belonging to $S$. In particular
\begin{align}
    w(\cal T) &= \sum_{T \in \cal T} w(T)\label{wei}\\
    w(\cal F^{x \to y}) &= \sum_{F \in \cal F^{x \to y}} w(F).
\end{align}

\subsection{Graphical representation}\label{grre}
We start by recalling the Kirchhoff formula (which is an application of the matrix-tree theorem).\\
Define $T_x^y$ to be the tree in $\cal T_x$ such that the edge between vertex $y$ and $y + 1/N$ is removed. It follows that $w(\cal T) = \sum_x w(\cal T_x)$ and $w(\cal T_x) = \sum_y w(T^y_x)$.\\
Using the matrix-tree theorem \cite{kir}, the stationary probability distribution $\rho^s$, $\rho^s L = 0$, appears as given by
\begin{equation}\label{Matrixtree}
\rho^s(x) = \frac{w(\cal T_x)}{\sum_{y}w(\cal T_y)} = \frac{\sum_y w(T_x^y)}{w(\cal T)}.
\end{equation}
We use it in Section \ref{stat} for taking the diffusion limit of the stationary distribution.\\

For now and next, we move from the single tree to the double tree level.
Recall Proposition \ref{p1} for the quasipotential \eqref{pp}.  For the random walk on the ring $\bbZ_N$, we have the following representation.
\begin{theorem}\label{TheFormula}
	When $\left\langle f\right\rangle^s =0 $,
\begin{equation}\label{gew}
V_f (x) =  \dfrac{ \sum_y  w\left(\cal F^{x\rightarrow y}\right)}{ w\left(\cal T\right)}f(y).
\end{equation}
\end{theorem}

\begin{proof}
    Define $g(x) =  \dfrac{ \sum_y  w\left(\cal F^{x\rightarrow y}\right)}{ w\left(\cal T\right)}f(y)$ for all $x$. 
    From \eqref{L}, 
    \begin{align*}
        Lg\,(x) = k(x, x+\frac{1}{N})&(g(x+\frac{1}{N}) - g(x)) + k(x, x-\frac{1}{N})(g(x-\frac{1}{N}) - g(x))
    \end{align*}
    and hence,
    \begin{align}
       w(\cal T) Lg\,(x) = k(x, x+\frac{1}{N}) \sum_y &\left(w(\cal F^{x+1/N \to y})- w(\cal F^{x \to y}) \right) f(y) + \label{bigsum}\\ 
       & k(x, x-\frac{1}{N}) \sum_y \left(w(\cal F^{x-1/N \to y}) - w(\cal F^{x \to y}) \right) f(y). \nonumber
    \end{align}
    We start with the case $y = x$ in the sum. We have that $\cal F^{x+1/N \to x} \subset \cal F^{x \to x}$. If $F \in \cal F^{x \to x} \setminus \cal F^{x+1/N \to x}$, then the edge $(x,x+1/N)$ does not belong to $F$. It follows that $k(x, x+1/N)w(F) = w(T)$ for some $T \in \cal T$. Remark that this tree is not rooted in $x$. In the same way we have that $\cal F^{x-1/N \to x} \subset \cal F^{x \to x}$ and if $F' \in \cal F^{x \to x} \setminus \cal F^{x-1/N \to x}$, then the edge $(x,x-1/N)$ does not belong to $F'$. Therefore, $k(x, x- 1/N)w(F) = w(T')$ for some $T' \in \cal T$. Remark again that $T'$ is not rooted in $x$. The other way around, for any $T \in \cal T$ that is not rooted in $x$, there exists $F \in \cal F^{x \to x}$ such that either $w(T) = k(x, x+1/N)w(F)$ or $w(T) = k(x, x-1/N)w(F)$. In this way we have a one-to-one correspondence between directed trees not rooted in $x$ and directed double rooted trees in the ring. As a consequence, 
    \begin{align*}
        k(x, x+\frac{1}{N}) \big(&w(\cal F^{x+1/N \to x}) - w(\cal F^{x \to x}) \big) f(x) \\
        &+ k(x, x-\frac{1}{N}) \left(w(\cal F^{x-1/N \to x}) - w(\cal F^{x \to x}) \right) f(x) = -\sum_{T \notin \cal T_x} w(T) f(x).
    \end{align*}
    Next, we consider what  remains in the sum \eqref{bigsum}, i.e.,
    \begin{align}
       k(x, x+\frac{1}{N}) \sum_{y \neq x} &\left(w(\cal F^{x+1/N \to y}) - w(\cal F^{x \to y}) \right) f(y) + \label{remaining} \\ 
       & k(x, x-\frac{1}{N}) \sum_{y \neq x} \left(w(\cal F^{x-1/N \to y}) - w(\cal F^{x \to y}) \right) f(y). \nonumber
    \end{align}
    Let $F \in \cal F^{x+1/N \to y}$ and $(x, x+1/N) \in F$. Then also $F \in \cal F^{x \to y}$.  Similarly, if $F \in \cal F^{x \to y}$ and $(x, x+1/N) \in F$, then $F \in \cal F^{x+1/N \to y}$.  The same holds for $x-1/N$ instead of $x+1/N$. \\
    Now let $F \in \cal F^{x+1/N \to y}$ and $(x, x+1/N) \notin F$. Suppose that $x$ is no root in $F$. Then there is a unique $F'\in \cal F^{x \to y}$ with $(x, x-1/N) \notin F'$ such that $k(x, x+1/N)w(F) = k(x, x-1/N)w(F')$. The same for $F \in \cal F^{x-1/N \to y}$ and $(x, x-1/N) \notin F$ such that $x$ is not a root. Then there is a unique $F'\in \cal F^{x \to y}$ with $(x,x+1/N) \notin F'$ such that $k(x, x-1/N)w(F) = k(x, x-1/N)w(F')$. We can also start with $F \in \cal F^{x \to y}$ with $(x,x+1/N) \notin F$, then there is $F' \in \cal F^{x-1/N \to y}$ without root in $x$ such that $(x,x-1/N) \notin F'$ and $k(x,x-1/N)w(F') = k(x, x+1/N)w(F)$. The same reasoning holds if we start with $F \in \cal F^{x \to y}$ with $(x,x-1/N) \notin F$. We have established the one--to--one correspondence. \\
    The only remaining terms in the sum \eqref{remaining} are
    \begin{align*}
        k(x,x+1/N)\sum_{y \neq x}\sum_{\substack{F \in \cal F^{x+1/N \to y}\\ (x, x+1/N) \notin F\\ x \text{ is root}}} w(F) + k(x,x-1/N)\sum_{y \neq x}\sum_{\substack{F \in \cal F^{x-1/N \to y}\\ (x, x-1/N) \notin F\\ x \text{ is root}}} w(F) .
    \end{align*}
    This is equal to $$ \sum_{\substack{T \in \cal T_y \\ y \neq x}} w(T)f(y).$$   
    Altogether we have that \eqref{bigsum} is equal to 
    \begin{align*}
        &-\sum_{T \notin \cal T_x}w(T)f(x) + \sum_{\substack{T \in \cal T_y \\ y \neq x}} w(T)f(y) \\ &= -\sum_{T \notin \cal T_x}w(T)f(x) +\sum_{\substack{T \in \cal T_y \\ y \neq x}} w(T)f(y) - \sum_{T \in \cal T_x} w(T)f(x) + \sum_{ T \in \cal T_x}w(T)f(x)\\
        & = -w(\cal T) f(x) +\sum_{T \in \cal T_y} w(T) f(y)\\
        &= -w(\cal T)f(x).
    \end{align*}
    The last equality follows from the assumption that $f$ has vanishing stationary average, and the Kirchhoff formula \eqref{Matrixtree}.\\
    We conclude that $w(\cal T) Lg\,(x) = -w(\cal T) f(x)$ and thus $Lg\,(x) = -f(x)$ for all $x$. Hence, $-(L_D^{-1} f)(x) = g(x) =  \dfrac{ \sum_y  w\left(\cal F^{x\rightarrow y}\right)}{ w\left(\cal T\right)}f(y)$.
\end{proof}

The proof of the above \eqref{gew} can be seen as an explicit proof of the matrix-forest theorem for the Markov jump process on $\bbZ_N$.
 For (more complicated) graphical expressions of the quasipotential of a general Markov jump process on a finite connected graph, see \cite{mathnernst}.

\subsection{Computation of quasipotential}\label{dop}
The graphical representation \eqref{gew} leads to an algorithm discussed in Appendix \ref{alg} for efficient computation of the quasipotential.  The backward generator $L$ takes the form
\begin{equation} \label{lr}
Lg\,(x) = k(x,x+\frac{1}{N})\,\left(g(x+\frac{1}{N}) - g(x)\right)+ k(x,x-\frac{1}{N})\,\left(g(x-\frac{1}{N}) - g(x)\right)
\end{equation}
on functions $g$ defined on $\mathbb{Z}_N$.\\ 
Let $u$ be an energy function defined on $[0,1]$ and write $\beta = \frac{1}{T}$  for the inverse of the absolute temperature $T$.
We use three types of transition rates:
\begin{equation} \label{rate1}
    k(x,x+\frac{1}{N}) = e^{\beta(u(x) - u(x+\frac{1}{N}))}e^{\frac{\varepsilon}{2}\frac{1}{N}},\quad
k(x,x-\frac{1}{N}) = e^{\beta(u(x) - u(x-\frac{1}{N}))}e^{-\frac{\varepsilon}{2}\frac{1}{N}}
\end{equation}
(unbounded rates with drift amplitude $\frac{\varepsilon}{N}$)
\begin{equation} \label{rate2}
    k(x,x+\frac{1}{N}) = e^{\frac{\beta}{2}(u(x) - u(x+\frac{1}{N}))}e^{\beta\frac{\varepsilon}{2}\frac{1}{N}} ,\quad
k(x,x-\frac{1}{N}) = e^{\frac{\beta}{2}(u(x) - u(x-\frac{1}{N}))}e^{-\beta\frac{\varepsilon}{2}\frac{1}{N}}
\end{equation}
(unbounded rates with drift amplitude $\beta\frac{\varepsilon}{N}$), and
\begin{equation} \label{rate3}
    k(x,x+\frac{1}{N}) = \frac{e^{\frac{\varepsilon}{2}\frac{1}{N}}}{1 + e^{-\beta(u(x)-u(x+\frac{1}{N}))}},\quad
k(x,x-\frac{1}{N}) = \frac{e^{-\frac{\varepsilon}{2}\frac{1}{N}}}{1 + e^{-\beta(u(x)-u(x-\frac{1}{N}))}}
\end{equation}
(bounded rates with drift amplitude $\frac{\varepsilon}{N}$).\\
Each time, when $\varepsilon=0$, the reversible probability distribution is the Gibbs distribution $\rho_\text{eq}(x) \propto \exp[-\beta u(x)]$. If $\varepsilon \neq 0$, there is a stationary current along the ring and we get a nonequilibrium steady state, which is not determined by the energy function $u$ solely.  If not specified, below, the $k(x,y)$ refer to any of the choices above.\\

\begin{figure}[h!]
	\RawFloats
	\begin{minipage}[t]{.4\textwidth}
		\centering
		\includegraphics[width=1.15\textwidth]{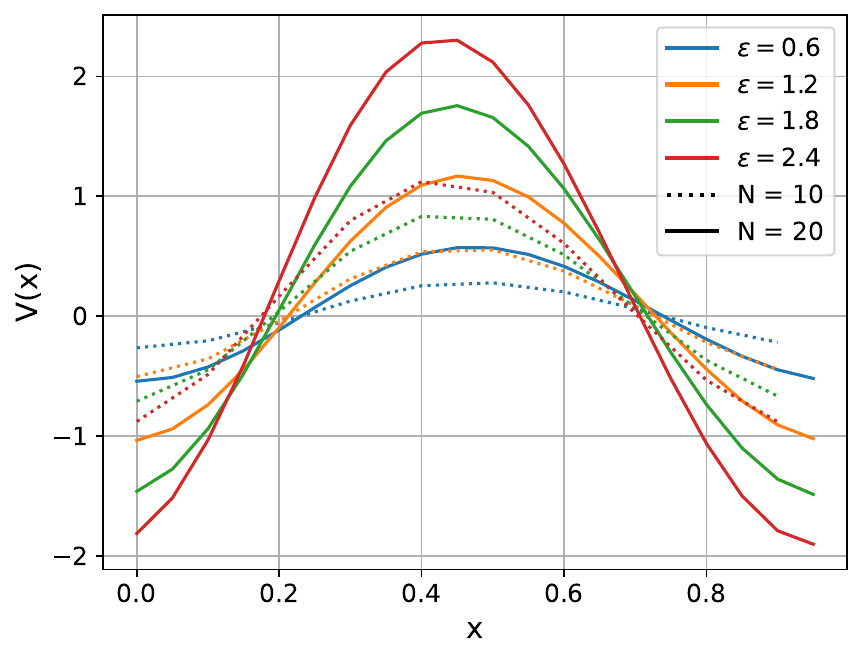}
		\caption{{\small Quasipotential as a function of $x$ at $\beta=1/2$ with rates \eqref{rate1}.}}\label{V(x) rate 1}
	\end{minipage}
	\hspace{1.5cm}
	\begin{minipage}[t]{.4\textwidth}
		\centering
		\includegraphics[width=1.15\textwidth]{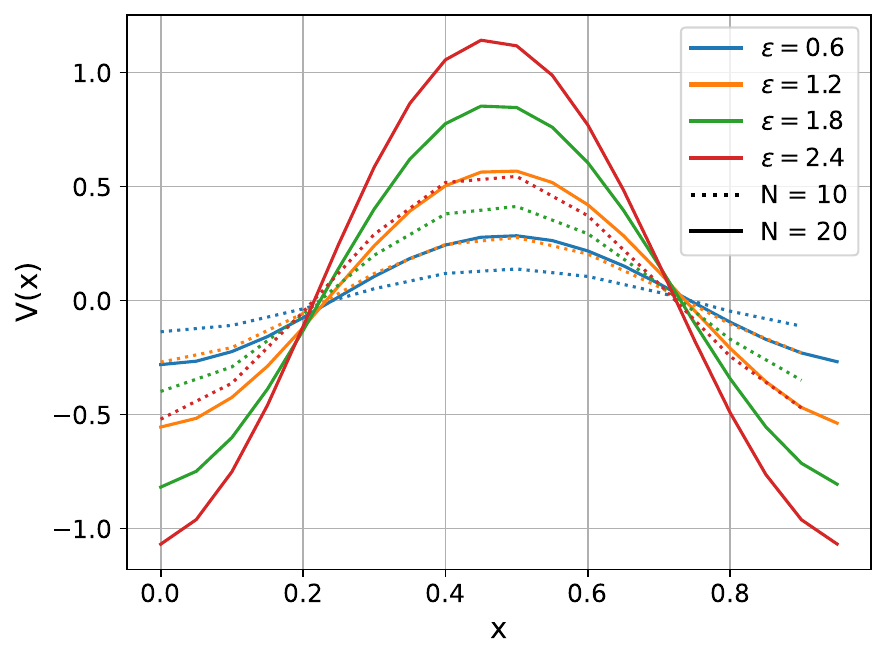}
		\caption{{\small Quasipotential as a function of $x$ at $\beta=1/2$ with rates \eqref{rate2}.}}\label{V(x) rate 2}
	\end{minipage}
	\hspace{0.5cm}
	\begin{minipage}[t]{.4\textwidth}
		\centering
		\includegraphics[width=1.15\textwidth]{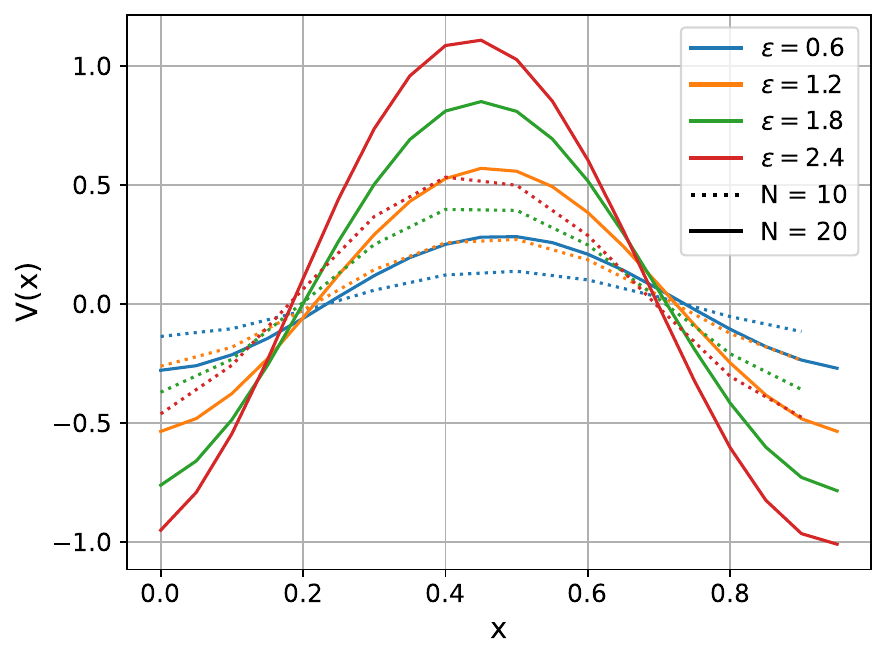}
		\caption{{\small Quasipotential as a function of $x$ at $\beta=1/2$ with rates \eqref{rate3}.}}\label{V(x) rate 3}
	\end{minipage} 
\end{figure}

The quasipotential $V = V_{q}$ is taken as in \eqref{pp} but for the special source $f=q$,
\begin{equation}\label{force}	
q(x) := h(x) - \langle h\rangle,\quad h(x) := -\ve \left[ k(x,x+\frac{1}{N}) -  k(x,x-\frac{1}{N})\right]
\end{equation}
which physically represents the Joule heating, i.e. the dissipated power.  
We show plots of $V(x),\, x\in [0,1)$ in Figs.~\ref{V(x) rate 1}--\ref{V(x) rate 2}--\ref{V(x) rate 3} for different $\varepsilon$ and $N$.  They correspond to the three choices of transition rates \eqref{rate1}--\eqref{rate2}--\eqref{rate3} with energy function $u(x) = 0.3 \sin(2 \pi x), \, x \in [0,1)$ and temperature $T = \frac{1}{\beta} = 2$.\\ 

It is remarkable  that the shape of the quasipotential $V$ is quite similar over different $N,\varepsilon$ and for different choices of transition rates. The function $V$ gets flatter for higher temperatures (not shown) and for smaller $\varepsilon$ and $N$.  The maximum is around $x=\frac{1}{2}$. The latter also determines the global shape of $V$.\\
Note that the curves with the same $\dfrac{\varepsilon}{N}$ almost coincide, which indicates the possibility of taking the continuum (diffusion) limit as indeed performed in Section \ref{dp}, and also explaining the similarity in shapes over the different curves.\\  We emphasize that the figures represent {\it exact} calculations; no approximations concerning temperature, size, driving, or energy  are needed.

\section{Heat capacities}\label{hc}
Following recent studies on thermal response for nonequilibrium systems \cite{pri,oon,kom2,eu,jir,noncal,sct}, we show here how to exactly compute the heat capacity 
\begin{equation}\label{heatc}
C(T) = \frac{\id}{\id T}\langle u \rangle - \left\langle \frac{\id V}{\id T} \right\rangle 
\end{equation}
where $V$ is the  quasipotential as in \eqref{pp}.  The point is that, from the graphical representation of its quasipotential $V$, we
obtain an efficient computation of the heat capacity $C(T)$ as well.\\

\begin{figure}[H]
	\RawFloats
\begin{minipage}[t]{.4\textwidth}
\includegraphics[width=1.25\textwidth]{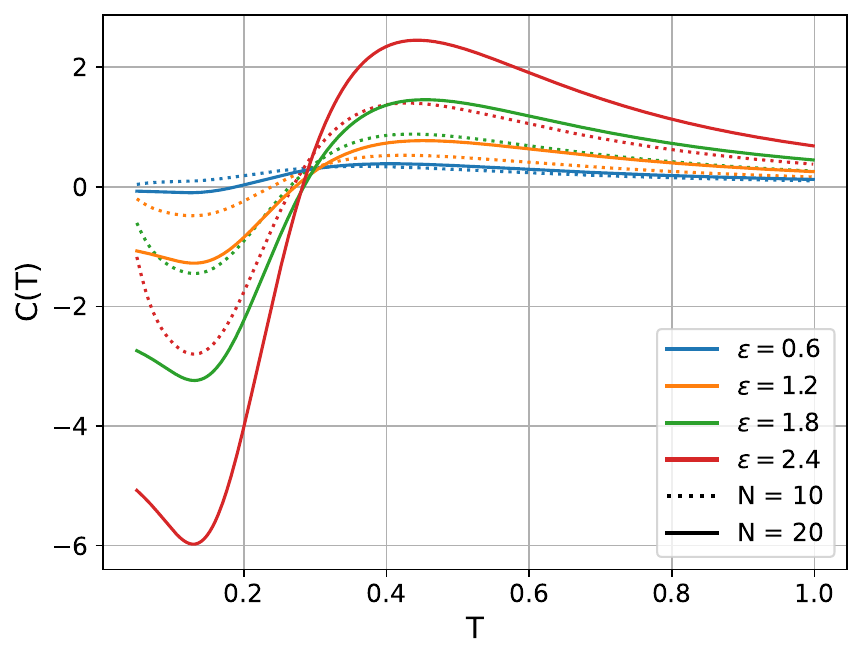}
\caption{{\small Heat capacity  using rate \eqref{rate1}.}}\label{C(T) rate 1 not close}
\end{minipage} 
	\hspace{1.5cm}
\begin{minipage}[t]{.4\textwidth}
		\includegraphics[width=1.2\textwidth]{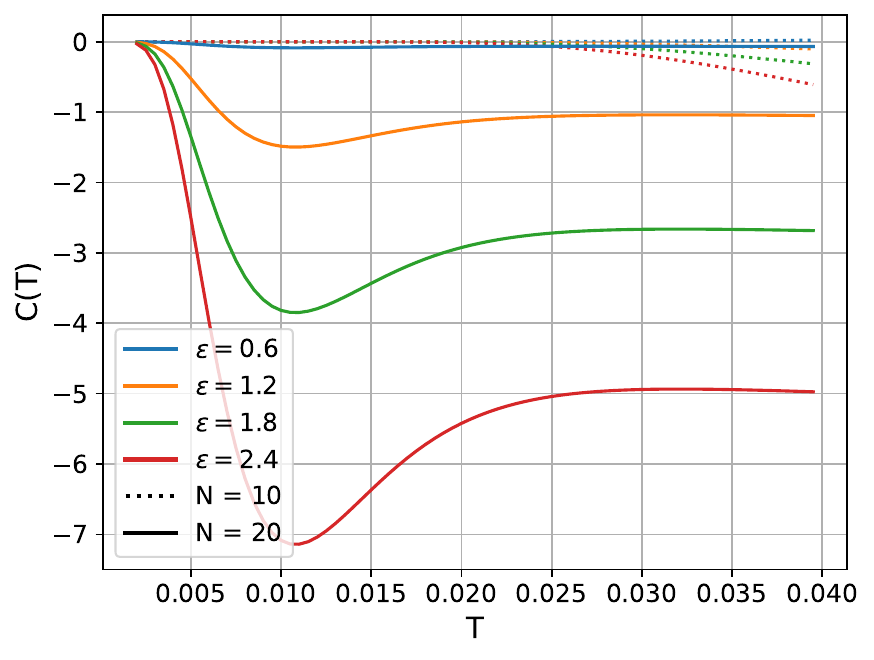}
		\caption{{\small Same for low temperatures.}}\label{C(T) rate 1 close to 0}
\end{minipage}
\end{figure}

In the figures below, we  use the energy function $u(x) = 0.3 \sin(2 \pi \frac{x}{N}), x\in \{0,1,\dots,N-1\}$ and we vary the ring size $N$ and the driving $\varepsilon$.\\
The heat capacity $C(T)$ using rate \eqref{rate1} is shown in Figs.~\ref{C(T) rate 1 not close} and \ref{C(T) rate 1 close to 0}.   Fig.~\ref{C(T) rate 1 close to 0} gives the heat capacity $C(T)$ for small temperatures $T$.    Interestingly, the heat capacity gets negative at low enough temperatures whenever $\varepsilon \neq 0$, which is a pure nonequilibrium effect.  It signifies, paradoxically, that the system needs to absorb net energy to reach lower energy states.  The extremes get bigger (in absolute value) when $\varepsilon$ increases.

\begin{figure}[H]
	\RawFloats
	\begin{minipage}[t]{.4\textwidth}
		\centering
		\includegraphics[width=1.2\textwidth]{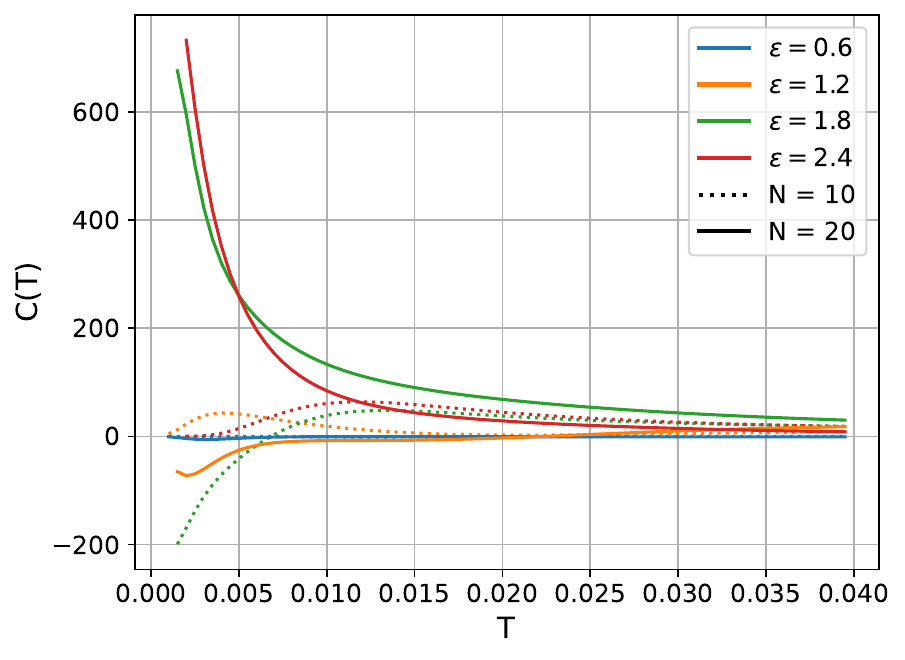}
		\caption{{\small Low temperature heat capacity using \eqref{rate2}. }}\label{C(T) rate 2 close to 0}
	\end{minipage}
	\hspace{1.5cm}
	\begin{minipage}[t]{.39\textwidth}
		\centering
		\includegraphics[width=1.2\textwidth]{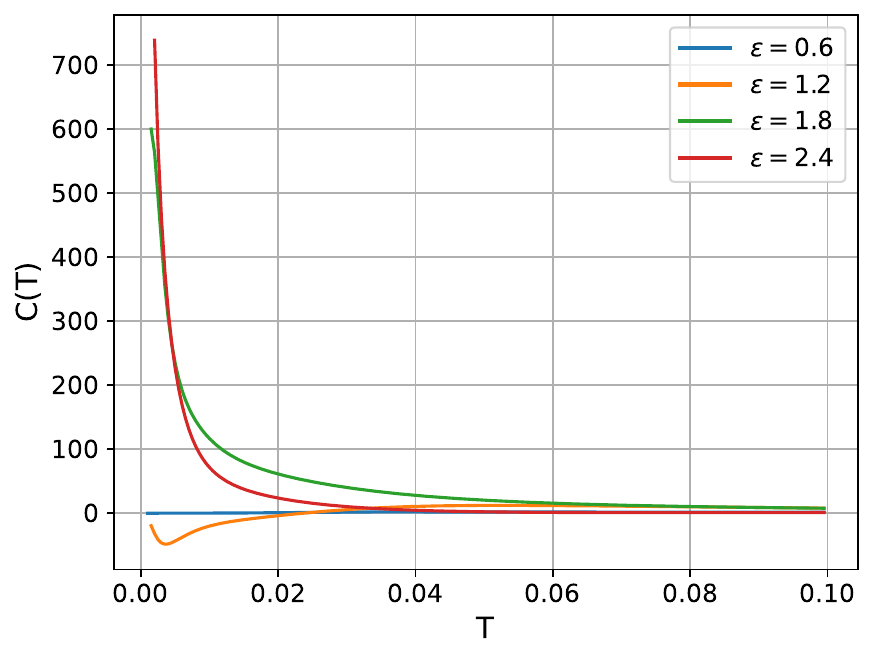}
		\caption{\small{Heat capacity with $N = 10 \ve$ using rate \eqref{rate2}.}}\label{C(T) N = 10eps rate 2}
	\end{minipage}   
\end{figure}

The heat capacity using rate \eqref{rate2} is given in Figs.~\ref{C(T) rate 2 close to 0} and \ref{C(T) N = 10eps rate 2}. The left shows $C(T)$ for temperatures $T$ close to zero.   Here, the heat capacity does not vanish at absolute zero (except in equilibrium $\varepsilon=0$) and appears monotone at larger $\varepsilon$.  To the right  is the heat capacity with fixed $N = 10 \varepsilon$, which is a good procedure for taking the continuum limit $N\uparrow\infty$.
 
\begin{figure}[H]
	\RawFloats
	\begin{minipage}[t]{.4\textwidth}
		\centering
		\includegraphics[width=1.2\textwidth]{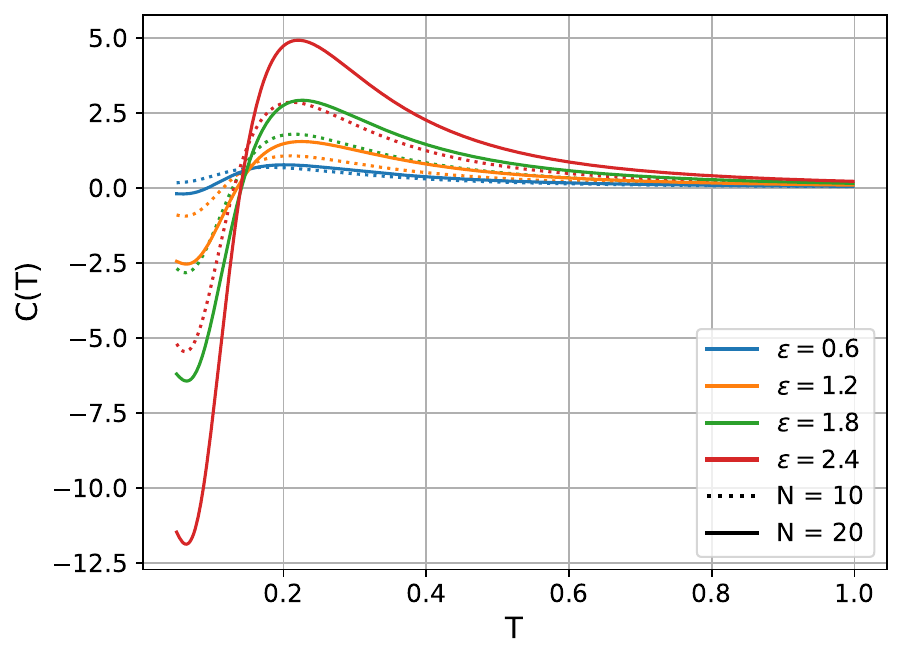}
		\caption{\small{Heat capacity using rate \eqref{rate3}}.}\label{C(T) rate 3 not close}
	\end{minipage} 
	\hspace{1.5cm}
	\begin{minipage}[t]{.4\textwidth}
		\centering
		\includegraphics[width=1.2\textwidth]{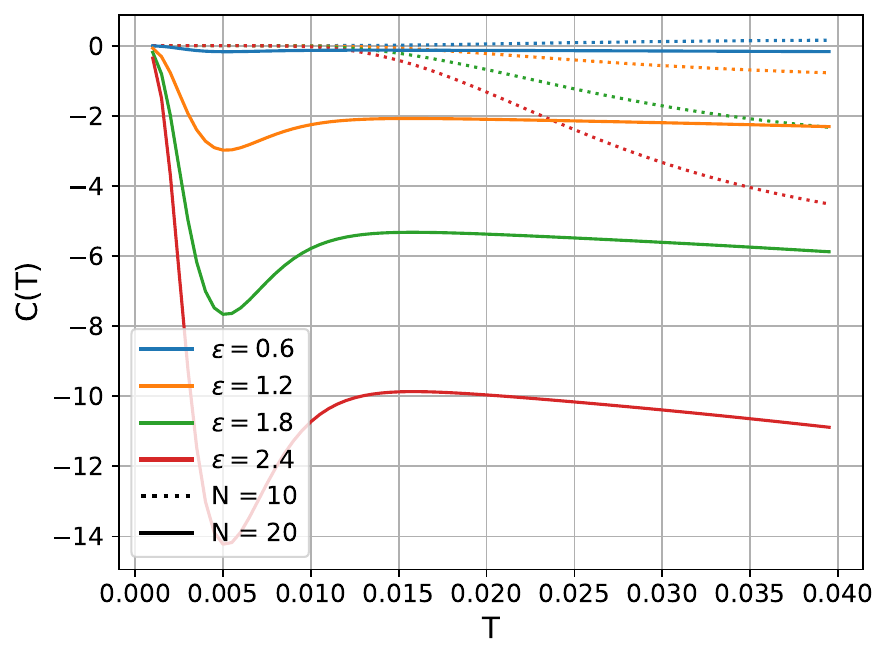}
		\caption{{\small Same for low temperatures.}}\label{C(T) rate 3 close to 0}
	\end{minipage}
\end{figure}

For the choice of rates \eqref{rate3}, the heat capacity is shown in Figs.~\ref{C(T) rate 3 not close}--\ref{C(T) rate 3 close to 0}. To the right, the heat capacity $C(T)$ is plotted at low temperatures $T$.  We see that the heat capacity always vanishes at absolute zero $T=0$.  The left plot in Fig.~\ref{C(T) rate 3 not close} is also very similar to Fig.~\ref{C(T) rate 1 not close}. 

 We emphasize that the heat capacity $C(T)\rightarrow 0$ as $T\downarrow 0$ in
Figs.~\ref{C(T) rate 1 close to 0} and \ref{C(T) rate 3 close to 0}, as illustrates a more
general extension of the Nernst heat theorem to nonequilibrium systems
\cite{ner,faz}. \\

 To the best of our knowledge, such precision in an exact
calculation of thermal response for a physically motivated nonequilibrium system close to
zero temperature has never been achieved,  \cite{code}.  The fact that our
computations reach so close to absolute zero shows their efficiency
and usefulness.\\
The behavior at high temperatures is easy to understand as the heat capacity needs to reach zero there. At any rate, the full parameter range can be explored with the method presented here, giving exact results for the heat capacity.  The only limiting factor is computer time which grows as  $N^4$.

\section{Diffusion limit from the Kirchhoff formula}\label{stat}
Consider the Markov diffusion process $x_t \in S^1$ (unit circle) defined via,
\begin{equation}\label{de}
\id x_t = [\ve - u'(x_t)] \,\id t + \sqrt{2T}\,\id W_t
\end{equation}
where $W_t$ is standard Brownian motion, with prefactor containing the temperature $T$.  It is again a driven stochastic model, with energy function $u$ and constant driving $\ve$.  It can be obtained  as the diffusion (continuous) limit of the random walk model on the ring that we have considered above.  We need to rescale time by $N^2$ while the mesh shrinks like $1/N$; see e.g.~\cite{davide}. More specifically, on smooth functions $g$ and for $x\in S^1$,
\begin{equation}\label{leg}
\beta^{-1}N^2\,Lg(\lfloor Nx\rfloor/N) = \bbL g(x),\quad \bbL g(x)= (-u'(x)+
\ve)g'(x) + \beta^{-1}g''(x)
\end{equation}
where $L$ is the backward generator \eqref{lr} of the Markov process on $\frac 1{N}\bbZ_N$ with  transition rates \ref{rate2}.
In that way, the stationary density $\rho$ of \eqref{de} can be obtained by considering the $N\uparrow\infty$ limit of the stationary distribution on the $N$ states of the (discrete) ring: formally, $\rho(x) = \lim_N N\rho^s(\lfloor Nx\rfloor/N),\, x\in S^1$.  The present section is devoted to that limit, repeating in somewhat simpler terms the calculation in \cite{davide}, starting from the Kirchhoff-matrix tree theorem.\\

Recall the tree notions from Section \ref{mft}, and the Kirchhoff formula \eqref{Matrixtree} in particular. 
That Kirchhoff formula  has been applied in \cite{heatb} to obtain the low-temperature asymptotics of $\rho^s$; other illustrations of that formula are contained in \cite{intro}.  We use \eqref{wei} with transition rates \eqref{rate2}.  For those rates the amplitude $\ve$ is divided by $N$, which is the correct diffusive scaling to reach \eqref{de}.\\

\begin{figure}[H]
    \centering
    \def\svgwidth{0.3\linewidth}        
    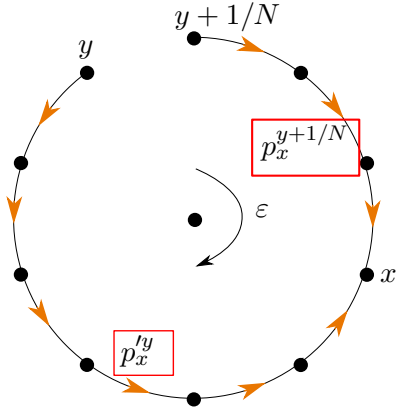
    \caption{{\small Different $T^y_x$ corresponding to $y$ and $x$ positions.}}\label{figurecase1-2-3}
\end{figure}

Remember from the beginning of Section \ref{grre} that $T_x^y$ denotes the tree in $\cal T_x$ such that the edge between vertex $y$ and $y + 1/N$ is removed.  Therefore, in $T^y_x$, there are paths from $y$ to $x$ and from $y+1/N$ to $x$. Hence, $T^y_x$ splits into  two parts: $p'^y_x$ denotes the path starting from $y$ with all edges directed (counter-clockwise) and ending in $x$, and we let $p^{y+1/N}_x$ be the path starting from $y+1/N$ where all edges are directed clockwise and ending in $x$; see Fig.~\ref{figurecase1-2-3}. 
From \eqref{Matrixtree}, 
\[w (\cal T_x)=\sum_y w(T_x^y)=\sum_y w(p'^y_x)w(p^{y+1/N}_x)\]
where $w(p'^y_x)$ is the product of all transitions rates of the edges belonging to the path from $y$ to $x$ and $p^{y+1/N}_x$ is the product of all transitions rates of the edges belonging to the path from $y + 1/N$ to $x$. 
If we use rates \eqref{rate2}, we get
\begin{align*}
    w (\cal T_x)=\sum_y e^{\frac{\beta }{2}(u(y)-u(x))}\,e^{-\frac{\beta}{2 N} \ve a}\, e^{\frac{\beta}{2}(u(y+\frac{1}{N})-u(x))}\,e^{\frac{\beta}{2 N} \ve b}.
\end{align*}
Here, for every $y$, $a+b=N-1$, where $a$ is the number of edges between $y$ and $x$, and  $b$ is the number of edges between $y+1/N$ and $x$.\\
In the limit of large $N$, 
\[\lim_{ N \rightarrow \infty}\frac{\beta}{2 N} \ve  a =\frac{\beta }{2 }\int_y^x\ve \, dz\]
where $x-y\geq 0$ and, similarly with $b=N-1-a$, 
\[
\lim_{ N \rightarrow \infty}\frac{\beta}{2 N} \ve  b=\frac{\beta \ve }{2 }-\frac{\beta }{2 }(x-y)\ve.
\]

Rescaling with $N$ to reach the stationary density, we conclude that
\begin{align}
\lim_{N \to \infty}	\frac 1{N} \sum_y w(T_x^y) =\int dy &e^{\beta(u(y) - u(x) +\varepsilon(y-x) +\frac{\varepsilon}{2}  )}.
\end{align}
That shows that the stationary distribution of \eqref{de} can be written as
\begin{equation}\label{dif}
    \rho(x)=\frac{1}{\cal Z}\int  \id y\, e^{\beta \, W(y,x)}
\end{equation}
where $W(y,x)=u(y)-u(x)+ (y-x)\ve +\frac{\ve}{2}$, 
and 
\[\cal Z=\int\, \int \id y \, \id x \, e^{\beta \,W(y,x)}\]
where the order on the circle is given by the direction of the driving with $\varepsilon$.\\
Another more direct derivation is contained in \cite{sdiv} where the stationary probability density for \eqref{de} is calculated by solving the Smoluchowski equation.  The above method can however be extended to an exact calculation for more general geometries.\\

\section{Diffusion limit of quasipotential}\label{dp}

The previous calculations can be extended to obtain the diffusion limit of the quasipotential \eqref{gew}.
Before doing that, we note that for the Markov diffusion process \eqref{de}, we define the quasipotential $\bbV_f$ as the smooth function on $S^1$ for which $\oint \id x\,\bbV_f(x) \rho(x)= 0$ (for probability density $\rho$ in \eqref{dif}) and $\bbL \bbV_f = -f$ (for backward generator $\bbL$ defined in \eqref{leg}), with source $f$ on $S^1$ having zero stationary expectation $\oint \id x\,f(x) \rho(x)= 0$ as well.  Therefore, that quasipotential $\bbV_f$ can be obtained from $R = -\bbV'$ with
\begin{equation}\label{rd}
[-u'(x)+\ve]\,R(x) + T R'(x) = f(x),\quad x\in S^1.
\end{equation}
For example, when $\ve=0$ (detailed balance), with $\oint\id y f(y) e^{-\beta u(y)} =0$, then the solution is
\begin{equation}\label{solr}
R(x) =  \beta e^{\beta u(x)}\int_0^x \id y f(y) e^{-\beta u(y)} - \beta \oint \id x e^{\beta u(x)}\int_0^x \id y f(y) e^{-\beta u(y)}.
\end{equation}
For $\ve\neq 0$, the solution $R$ can  be obtained in the same format as \eqref{dif}; see Eq.~3 in \cite{sdiv}.\\
In principle, the graphical representation \eqref{gew} gives an algorithm to approximate those solutions more generally.   We will not actually carry that through by hand, while the Figs.~\ref{V(x) rate 2}--\ref{V(x) rate 3} of course give such an approximation for specific source functions (as in \eqref{force}).  Coming from the quasipotential on the ring, one needs to be careful in rescaling $f$ with $N^2$, as already clear from \eqref{vv}, which reads $V''(x) = -N^2 f$ for large $N$.\\
We limit ourselves here to presenting the approximating formul{\ae}.\\

\begin{figure}[H]
    \centering
 \def\svgwidth{0.6\linewidth}        
    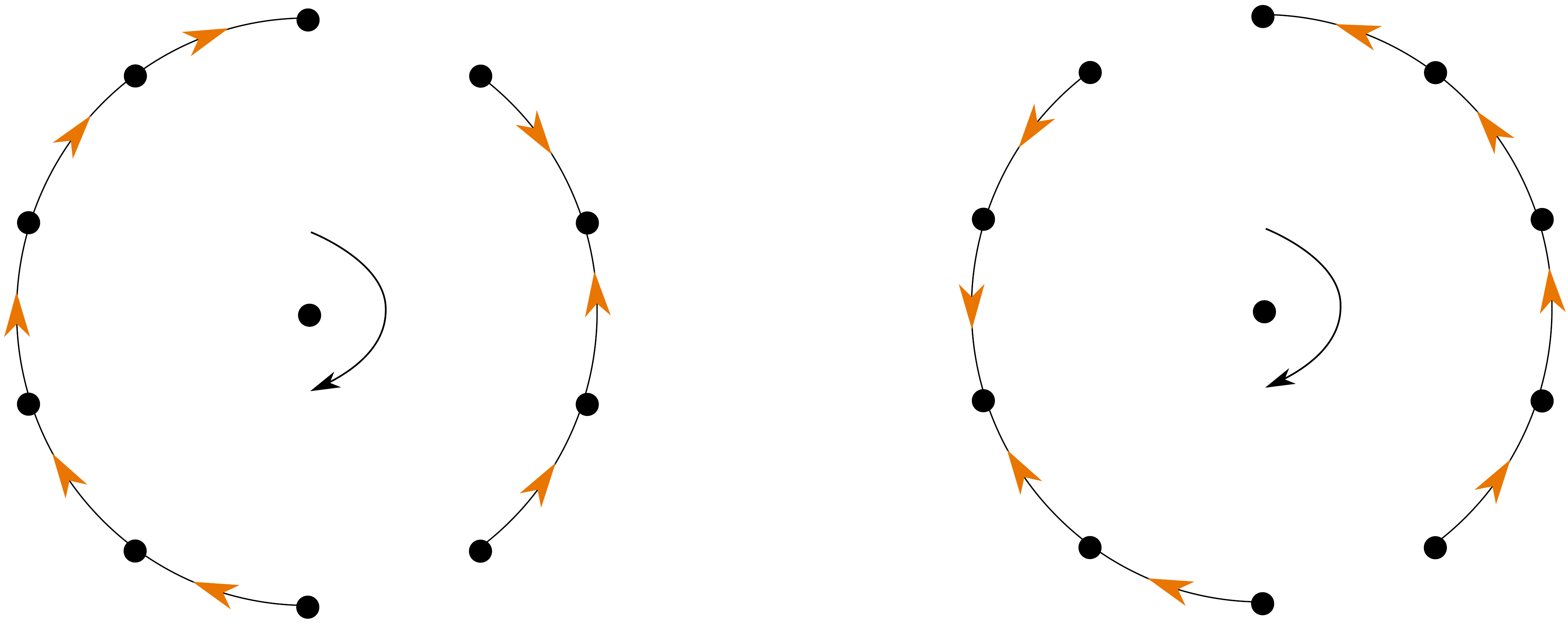
    \caption{{\small An element of the set $\cal F^{x\to x}\setminus \cal F^{x+1/N\to x}$ on the left side and an element of the set $\cal F^{x\to x}\setminus \cal F^{x-1/N\to x}$ on the right side.}}\label{y=x}
\end{figure}

To start,  consider 
\[V_+(x)= \bbV(x+\frac{1}{N})-\bbV(x), \qquad V_-(x)= \bbV(x)-\bbV(x-\frac{1}{N}).\]
From the graphical expression \eqref{gew}, we  have
\begin{align}\label{V+V-}
V_+(x)&=\frac{1}{w(\cal T)}\sum_y\big(w(\cal F^{x+1/N\to y})-w(\cal F^{x\to y})\big)f(y)\notag\\
V_-(x)&=\frac{1}{w(\cal T)}\sum_y\big(w(\cal F^{x\to y})-w(\cal F^{x-1/N\to y})\big)f(y).
\end{align}
First, take the term $y=x$ in $V_+(x)$.  $\cal F^{x\to x}$ is the set of all double trees where $x$ is the root of one part and the second part is rooted in some $r$.  $\cal F^{x+1/N\to x}$ is the set of all   double trees where $x$ and $x+1/N$ are in the same part and $x$ is the root. Hence, $\cal F^{x+1/N\to x}\subset \cal F^{x\to x}$ and  
\begin{align*}\label{V+}
\big(w(\cal F^{x+1/N\to x})-w(\cal F^{x\to x})\big)f(x)&= - w(\cal F^{x\to x}\setminus \cal F^{x+1/N\to x})f(x)
\end{align*}
where $\cal F^{x\to x}\setminus \cal F^{x+1/N\to x}$ is the set of all double trees  on the ring with the edge $\{x, x+\frac{1}{N}\}$ being removed and with a second gap $\{z, z+\frac{1}{N}\}$ for some $z$. The same holds for $V_-(x)$; see \fig \ref{y=x}. 

Referring to \fig\ref{y=x}, the double trees are either rooted in $r$ or in $x$. We introduce a notation for the two parts in the double tree. In a double tree with a gap between $z$ and $z+1/N$ we use the notation $\tau^z_{r}$ and $\tau^z_{r'}$ to denote the two parts, one rooted in $r$ and the other part rooted in $r'$. We have
\begin{align}
-w(\cal F^{x\to x}\setminus \cal F^{x+1/N\to x})&=-\sum_{ z=x+\frac{1}{N}}^{x-\frac{1}{N}} w(\tau^z_x)\,  \sum_{r=x+\frac{1}{N}}^{z} w(\tau^z_r)\notag\\
w(\cal F^{x\to x}\setminus \cal F^{x-1/N\to x})&=\sum_{z= x}^{x-\frac{2}{N}} w(\tau^z_x)\,  \sum_{r=z+\frac{1}{N}}^{ x- \frac{1}{N}} w(\tau^z_r).
\end{align}

\begin{figure}[H]
    \centering
    \def\svgwidth{0.8\linewidth}        
    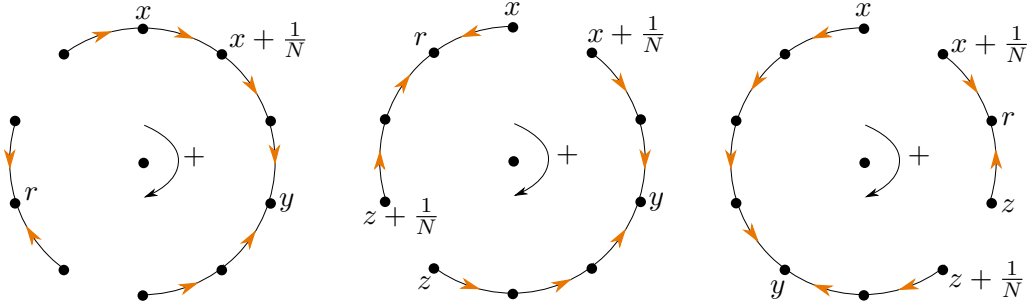
    \caption{{\small On the left side there is a double tree from the set $\cal F^{x\to y} \cap \cal F^{x+1/N\to y}$. In the middle there is a double tree where one of the gaps is between $x$ and $x+1/N$ and where $x+1/N$ is  located in the same tree as $y$. On the right side there is a double tree where one of the gaps is between $x$ and $x+1/N$ and $x$ is  located in the same tree as $y$.}}\label{v'}
\end{figure}

Next, we take  the terms with $y\not = x$ in $V_+(x)$.  $\cal F^{x\to y}$ is the set of double trees where $y$ and $x$ are located in the same part and $y$ is the root, while the second part is rooted in some $r$. These double trees are created by two gaps in the ring. The intersection of $\cal F^{x\to y}$ and $\cal F^{x+1/N\to y}$ is the set of  double trees including the edge $\{x,x+1/N\}$.  Hence, the double trees that are not in this intersection have a gap between $x$ and $x+1/N$, see Fig.~\ref{v'}. 

\begin{figure}[H]
    \centering
      \def\svgwidth{0.8\linewidth}        
    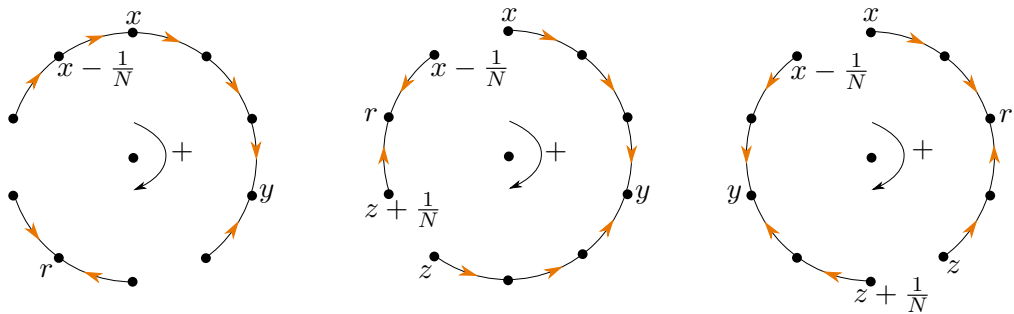
    \caption{{\small On the left side there is a double tree in $\cal F^{x-1/N\to y} \cap \cal F^{x\to y}$. In the middle, there is a double tree where one of the gaps is in between $x-1/N$ and $x$, and $x$ is located in the same tree as $y$. On the right side, there is a double tree where one of the gaps is between $x-1/N$ and $x$, and $x-1/N$ is located in the same tree as $y$.}}\label{v-}
\end{figure}

Suppose the second gap is between $z$ and $z+\frac{1}{N}$. Fix $y\not=x$, 
\begin{align}\label{diffv+}
 w(\cal F^{x+1/N\to y})- w(\cal F^{x\to y})= \sum_{ z=y}^{x-\frac{1}{N}}w(\tau^z_y)\, \sum_{r=z+\frac{1}{N}}^{x}w(\tau^z_r)
  -\sum_{z=x+\frac{1}{N}}^{y-\frac{1}{N}}w(\tau^z_y) \, \sum_{r= x+\frac{1}{N}} ^{ z}w(\tau^z_r).
 \end{align}
 In the same way when the gap is between  $x-\frac{1}{N}$ and $x$, see \fig \ref{v-}.
\begin{align}\label{diffv-}
 w(\cal F^{x\to y})- w(\cal F^{x-\frac{1}{N}\to y})= \sum_{ z=y}^{ x-\frac{2}{N}}w(\tau^z_y)\, \sum_{r=z+\frac{1}{N}}^{ x-\frac{1}{N} }w(\tau^z_r)
  -\sum_{z=x}^{ y-\frac{1}{N}}w(\tau^z_y) \, \sum_{r=x}^{ z}w(\tau^z_r).
 \end{align}
Hence, 
\begin{align}\label{tV+V-}
V_+(x)&=\frac{1}{w(\cal T)}\sum_{y\not=x}\bigg(\sum_{ z=y}^{x-\frac{1}{N}}w(\tau^z_y)\, \sum_{r=z+\frac{1}{N}}^{x}w(\tau^z_r)
  -\sum_{z=x+\frac{1}{N}}^{y-\frac{1}{N}}w(\tau^z_y) \, \sum_{r= x+\frac{1}{N}} ^{ z}w(\tau^z_r)\bigg)f(y)\notag\\
  &\quad -\frac{1}{w(\cal T)}\,f(x)\, \sum_{z=x+\frac{1}{N}}^{x-\frac{1}{N}}w(\tau^z_y) \, \sum_{r= x+\frac{1}{N}} ^{ z}w(\tau^z_r)\notag\\
V_-(x)&=\frac{1}{w(\cal T)}\sum_{y\not=x}\bigg(\sum_{ z=y}^{ x-\frac{2}{N}}w(\tau^z_y)\, \sum_{r=z+\frac{1}{N}}^{ x-\frac{1}{N} }w(\tau^z_r)
  -\sum_{z=x}^{ y-\frac{1}{N}}w(\tau^z_y) \, \sum_{r=x}^{ z}w(\tau^z_r)\bigg)f(y)\notag\\
  &\quad + \frac{1}{w(\cal T)}\,f(x)\,\sum_{ z=x}^{ x-\frac{2}{N}}w(\tau^z_x)\, \sum_{r=z+\frac{1}{N}}^{ x-\frac{1}{N} }w(\tau^z_r).
\end{align}

 The above formul{\ae} can be evaluated numerically, and  yield integrals after proper rescaling with $N\uparrow\infty$, to reach an approximation to the solution of  \eqref{rd}. We do not write down the resulting integrals. We believe however that the graphical representation of the solution of Poisson equations \eqref{pp} can be a more general tool to prove boundedness; see also \cite{mathnernst, poiss}.  Instead, we limit ourselves (in the next section) to show how that approach is helpful to obtain the boundedness of the quasipotential.

 \section{Boundedness of quasipotential}\label{nera}
 The previous formul{\ae} \eqref{tV+V-} and graphical representations allow evaluating the boundedness of the quasipotential. It suffices in fact to give  a bound for the difference of  quasipotentials over any edge.  We treat here two (extreme) cases of transition rates. \\

 We take the case of large $\ve>0$ in e.g. \eqref{rate3}, so that for simplicity we put $u(x)=0$ at all vertices $x$, while 
 \begin{equation}\label{rate3u0}
 k(x,x+\frac{1}{N})=e^{\frac{\ve}{2}\frac{1}{N}}, \qquad k(x,x-\frac{1}{N})=e^{-\frac{\ve}{2}\frac{1}{N}}.
 \end{equation}
Then, obviously, the maximal weight of  rooted spanning trees is obtained when all edges are directed clockwise (in the same direction as the bias $\ve$). There are  $N$ such rooted spanning  trees, each with weight $e^{\frac{\ve}{2}\frac{N-1}{N}}$. On the other hand, the maximum weight of double--rooted trees (with $N-2$ edges) is $e^{\frac{\ve}{2}\frac{N-2}{N}}$. Therefore, for every $\ve$, $V_+$ and $V_-$ in \eqref{tV+V-} are bounded,
\begin{align*}
V_+(x), \, V_-(x) \leq\sum_y \frac{N(N-1)e^{-\frac{\ve}{2}\frac{N-2}{N}}}{N  e^{-\frac{\ve}{2}\frac{N-1}{N}}} f(y)\leq \,N^2  \, e^{-\frac{\ve}{2}\frac{1}{N}} \,\, f_\text{max},\quad f_\text{max} := \max_y|f(y)|. 
\end{align*}
The $N^2$ is often too big as a prefactor (as we did not take into account the centrality of $f$).  For example, suppose that $f(x)=0$ whenever $x\not=x^*, x^*+1$, and  $\Delta f:=|f(x^*+1)-f(x^*)|$.  Then, 
 \[V_+(x),  \, V_-(x)\leq \frac{N(N-1)e^{-\frac{\ve}{2}\frac{N-2}{N}}}{N  e^{-\frac{\ve}{2}\frac{N-1}{N}}} \Delta f\leq N  e^{-\frac{\ve}{2}\frac{1}{N}}\Delta f\]
and if  for every $x$, $k(x,x+\frac{1}{N})=k(x,x-\frac{1}{N})=1$, then
 \begin{align*}
V_+(x),  \, V_-(x)\leq \frac{N^2 \Delta f}{N^2}\leq \Delta f. \end{align*}

For going to the other extreme (looking back at the rates \eqref{rate3} again), take $\ve =0$ and put for every $x\not=x^*$, $u(x)=0$ while $u(x^*)=u$.  That is 
 \[k(x^*,x^*+\frac{1}{N})= k(x^*,x^*-\frac{1}{N})=\frac{1}{1+e^{-\beta u}}, \qquad k(x,x+\frac{1}{N})=k(x,x-\frac{1}{N})=1.\]
 In that case $w(\cal T)=\frac{1}{1+e^{-\beta u}}(N(N-1)+2)+(N-2)$ and  
\begin{align*}
V_+(x^*), \, V_-(x^*) \leq \frac{N(N-1)}{\frac{1}{1+e^{-\beta u}}(N(N-1)+2)+(N-2)}\,N  f_{\text{max}}\leq N (1+e^{-\beta u})\, \,f_{\text{max}}
\end{align*}
and
 \begin{align*}
V_+(x), \, V_-(x) \leq \frac{(N^2-1)\frac{1}{1+e^{-\beta u}}+1}{\frac{1}{1+e^{-\beta u}}(N(N-1)+2)+(N-2)}\,N  f_{\text{max}}\leq N \frac{1}{1+e^{-\beta u}}\, \,f_{\text{max}}.
\end{align*}

\section{Conclusions}\label{con}

The computation of {\it excesses}, the temporally--accumulated difference
\begin{equation}\label{acc}
\int_0^\infty\id t\,\left[ \langle g(X_t)\,|\,X_0=x\rangle - \langle
g\rangle \right]
\end{equation}
between finite--time expectations  $ \langle g(X_t)\,|\,X_0=x\rangle$
and their  limit $\langle g\rangle  = \lim_{t\uparrow \infty}\langle
g(X_t)\,|\,X_0=x\rangle$,  is important for statistical mechanical
applications of stochastic processes $X_t, t\geq 0$.  In the present
paper, for $X_t$ we have considered a driven random walker on a ring with $N$
vertices, and we obtained a graphical representation (in terms of
combining weights of two trees) leading to efficient computation of
such a difference \eqref{acc}, called the quasipotential.  We have related that
quasipotential to the Drazin--inverse of the backward generator. \\
As an application to thermal physics and by choosing the appropriate
function $g$, we have computed the heat capacity of an ideal gas of
driven walkers.  We have demonstrated that it is perfectly possible,
thanks to the developed algorithm, to study its low-temperature behavior
(satisfying an extended Nernst heat theorem in the case of bounded
transition rates) and its dependence on driving and on the size of the
ring.  We have also shown that a continuum limit is possible,
giving an overdamped diffusion for which the quasipotential can be
given explicitly.\\
We expect that the combination of graph theory with nonequilibrium
statistical mechanics of Markov jump processes, as tried here, will give
rise to further interesting applications, in particular by extending
this work to metric graphs and to active (or persistent) random walkers.\\

\noindent{\bf Acknowledgment}: We are grateful to  Prof.~Christian Maes for suggesting the problem, many discussions and help in editing the paper. We thank  Mathias Smits and Simon Krekels for substantial help with the numerical work.

\newpage
\appendix

\section{Algorithm}\label{alg}

The applied algorithm \cite{code} finds  all possible rooted trees and all possible double rooted trees in $\cal F^{x\rightarrow y}$ for the graph $\mathbb{Z}_N$. Below we add explanations to make the code more readable. 

The algorithm presents the vertices as numbers $0$ up to $N-1$. It displays a rooted tree or double rooted tree by an array $A$ with $N$ entries which can be either $0$, $1$ or $-1$. If the $i$-th spot of the array is $0$ (so if $A[i-1] = 0$), there is no edge between the vertices $i-1$ and $i$. If $A[i-1] = 1$, there is a  edge $(i-1,i)$. And if $A[i-1] = -1$, there is a directed edge $(i,i-1)$. For example, the following four arrays
\begin{align*}
[1,0,0,1], \quad [1,-1,0,0], \quad [1,0,1,0], \quad [1,0,-1,0]
\end{align*}

\begin{figure}[H]
    \centering
    \def\svgwidth{0.8\linewidth}
    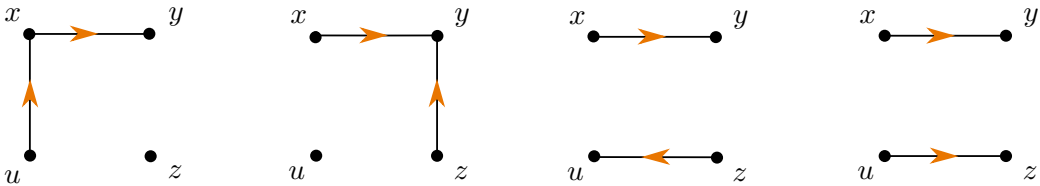
    \caption{{\small There are four forests in the set $\cal F^{x\rightarrow y}$. }}\label{forests}
\end{figure}

represent the double trees from ring $\bbZ_4$ where $x=0,\, y=1, \, z=2 $ and $u=3$ (see Fig.~\ref{forests}). Remark that there are always two zeroes in the arrays. This is because in $\mathbb{Z}_N$ a double tree in $\cal F^{x\rightarrow y}$ has $N-2$ edges. 
Since a tree has $N-1$ edges, only one entry in the array presenting a tree is zero. The other spots can be either $1$ or $-1$. 

\section{Laplacian matrix and its index}\label{lap}
The backward generator $L$ is a specific case of a Laplacian matrix when the underlining graph is weighted, directed, and does not have self-edges. In \cite{che2006}, the authors study this class of Laplacian matrices and in their Proposition 12, it is shown  that the index of the Laplacian matrix is one.  For more, see also \cite{adi}. Here, we are proving it for our special case of the  random walker on the ring.\\

\begin{lemma} \label{index}
	The index of the backward generator $L$ on the ring, see \eqref{lr}, is always one: ${\text ind}(L)=1$.
\end{lemma}
\begin{proof}
	 We use the notation $L_{ij}$ to denote the matrix element of $L$ on row $i$ and column $j$ where $i,j\,= 1, 2,..., N$. The matrix element $L_{ij}$ is zero for 
	\[\begin{cases}
	3 \leq j \leq N-1&\text{ if } i = 1\\
	1 \leq j < i-1\text{ or } i+1 < j \leq N&\text{ if } i \neq 1\text{ and }i \neq N\\
	2 \leq j \leq N-2&\text{ if } i = N.
	\end{cases}\]
	Since we use nonzero transition rates, every element of the form $L_{i,i-1}$ or $L_{i,i+1}$ is certainly non-zero. We also know that the sum over the elements of a row is zero.
	Let $M$ be the squared submatrix of $L$ where the first row and column are removed. The matrix $M$ is a $(N-1)\times(N-1)$-matrix having elements $M_{ij}$.
	Remark that in every $j$-th column of $M$, the element $M_{j+1,j}$ is nonzero and $M_{ij}$ with $i > j+1$ is zero. \\
	Proving that $M$ has full rank implies that $L$ has rank $N-1$. We do this by showing that the determinant of $M$ is nonzero. 
	We apply row operations on the matrix $M$. The fact that the determinant of $M$ is zero or nonzero will not change by doing that. We use the notation $M^{(i)} \xrightarrow{R_j \leftrightarrow R_k} M^{(i+1)}$ for exchanging the $j$-th row and the $k$-th row of $M^{(i)}$; the resulting matrix is called $M^{(i+1)}$. The notation $M^{(i)} \xrightarrow{R_j \to R_j/(M^{(i)})_{kl}} M^{(i+1)}$ indicates that we divide the $j$-th row of $M^{(i)}$ by $(M^{(i)})_{kl}$ and we call the resulting matrix $M^{(i+1)}$.  Finally, with $M^{(i)} \xrightarrow{R_{j} \to R_{j} - \lambda R_k } M^{(i+1)}$ we mean to subtract $\lambda$ times the $k$-th row of $M^{(i)}$ from the $j$-th row of $M^{(i)}$ and call the resulting matrix $M^{(i+1)}$. We apply the following row operations on the matrix $M$:
	\begin{align*}
	M &\xrightarrow{R_1 \leftrightarrow R_2} M^{(1)} \xrightarrow{R_1 \to R_1/(M^{(1)})_{11}} M^{(2)} \xrightarrow{R_2 \to R_2 - (M^{(2)})_{21} R_1 } M^{(3)} \\
	& \xrightarrow{R_2 \leftrightarrow R_3} M^{(4)} \xrightarrow{R_2 \to R_2/(M^{(4)})_{22}} M^{(5)} \xrightarrow{R_3 \to R_3 - (M^{(5)})_{32} R_2 } M^{(6)} \\
	& \cdots \rightarrow \cdots \\
	&\xrightarrow{R_i \leftrightarrow R_{i+1}} M^{(3i-2)} \xrightarrow{R_i \to R_i/(M^{(3i-2)})_{ii}} M^{(3i-1)} \xrightarrow{R_{i+1} \to R_{i+1} - (M^{(3i-1)})_{i+1,i} R_i } M^{(3i)}\\
	&\cdots \rightarrow \cdots.
	\end{align*}
	Remark that for every $i$, $(M^{(3i)})_{ii} = 1$ and $(M^{(3i)})_{ji} = 0$ for $j >i$. We end up with an upper square diagonal matrix.   The diagonal has no zeros. That shows that $M$ has full rank. As a consequence, $L$ has rank $N-1$.\\
	Next, we show that $L^2$ has rank $N-1$ as well. We have $\text{rank}(L^2) \leq \text{rank}(L)$. The matrix $L^2$ has elements $(L^2)_{ij}$ which are zero for
	\[\begin{cases}
	4 \leq j \leq N-1&\text{ if } i = 1\\
	5 \leq j \leq N&\text{ if } i = 2\\
	1 \leq j < i-2\text{ or } i+2 < j \leq N&\text{ if } i \notin \{1,2,N-1,N\}\\
	2 \leq j \leq N-4&\text{ if } i = N-1\\
	3 \leq j \leq N-3&\text{ if } i = N.
	\end{cases}\]
	Define $M$ now as the square $(N-1)\times(N-1)$ submatrix of $L^2$, where we remove the first column and the last row.
	Remark that in every $j$-th column of $M$ with $j<N-3$, the element $M_{j+2,j}$ is nonzero and $M_{ij}$ with $i \geq j+3$ is zero. Therefore, in each column, we can move the non-zero element up, as we did to prove that $L$ has rank $N-1$, and use that element to make the elements under the diagonal zero. Remark that while doing this, we do not change the value of $M_{j+2,j}$ in the other columns. Then we continue in the second column and we do the same process there and also in the other columns. In this way we end up with an upper--diagonal matrix with on the diagonal only the number $1$. The determinant of this matrix is nonzero and thus $M$ has full rank. That proves that $L^2$ has rank $N-1$.\\
	We conclude that $\text{ind}(L)  = 1$. 
\end{proof}

\end{document}